\documentclass[twoside, 12pt, a4paper]{article}
\usepackage{enumitem}
\usepackage{relsize}
\usepackage{eucal}
\usepackage{url}
\usepackage{fullpage}
\usepackage[perpage]{footmisc}
\makeatletter\def\@seccntformat#1{\protect\sc{\@ifundefined{#1@cntformat}{\csname the#1\endcsname\quad}%
{\csname #1@cntformat\endcsname}}}
\def\section@cntformat{\S\thesection.\ }
\def\subsection@cntformat{\thesubsection.\ }
\makeatother

\makeatletter
\renewcommand\section{\@startsection {section}{1}{\z@}%
                                   {-3.5ex \@plus -1ex \@minus -.2ex}%
                                   {2.3ex \@plus.2ex}%
                                   {\centering\normalfont\large\scshape}}

\makeatother

\makeatletter
\renewcommand\subsection{\@startsection {subsection}{1}{\z@}%
                                   {-3.5ex \@plus -1ex \@minus -.2ex}%
                                   {2.3ex \@plus.2ex}%
                                   {\normalfont\large}}

\makeatother
\usepackage{amsmath}
\usepackage[lite, numeric,sorted,traditional-quotes,compressed-cites]{amsrefs}
\usepackage{amssymb}
\usepackage{mathrsfs}
\usepackage[all]{xy}
\usepackage{bbm}
\usepackage[amsmath, thmmarks]{ntheorem}

\newfont{\sehrgrossen}{cmfib8 scaled 4300}

\newcommand{\mor}[1]{{\text{Hom}}_{#1}}

\newcommand{\id}[1]{{\text{id}}_{#1}}
\newcommand{\Mgm}[1]{M_{\text{gm}}(#1)}
\newcommand{\dime}[2]{\text{dim}_{#1}(#2)}
\newcommand{\class}[1]{\text{cl}_{#1}}
\newcommand{\clas}[1]{\text{l}_{#1}}
\newcommand{\Alt}[2]{\text{Alt}^{#1}(#2)}
\newcommand{\Sym}[2]{\text{Sym}^{#1}(#2)}
{\theorembodyfont{\sl} \theoremheaderfont{\scshape} \theoremseparator{.}
\newtheorem{theorem}{\textsc{Theorem}}[section]

\newtheorem{cor}[theorem]{Corollary}
\newtheorem{proposition}[theorem]{Proposition}
\newtheorem{lemma}[theorem]{Lemma}

}
{\theorembodyfont{\sl} \theoremheaderfont{\fontshape{sc}\selectfont} \theoremseparator{.}
\newtheorem{app-theorem}{\textsc{Theorem}}[section]

\newtheorem{app-cor}[app-theorem]{Corollary}
\newtheorem{app-proposition}[app-theorem]{Proposition}
\newtheorem{app-lemma}[app-theorem]{Lemma}
\newtheorem*{app-theoremn}{\underline{Theorem}}
\newtheorem*{app-theoremmain}{\underline{The Main Theorem}}
\newtheorem*{no-theorem}{Theorem}
}
{\theorembodyfont{\rmfamily} \theoremheaderfont{\scshape} \theoremseparator{.}
\newtheorem{remark}[theorem]{Remark}

\newtheorem{example}[theorem]{Example}

\newtheorem{definition}[theorem]{Definition}

\newtheorem{notations}[theorem]{Notations}

}

{\theoremstyle{nonumberplain}
\theoremheaderfont{\scshape}
\theorembodyfont{\normalfont}
\theoremseparator{.}
\theoremsymbol{\ensuremath{_\proofSymbol}}
\RequirePackage{amssymb}

\qedsymbol{\ensuremath{_\blacksquare}}
\theoremclass{LaTeX}}

\newenvironment{proof}{\textsc{Proof}.}{\hspace*{\fill}
$\square$\\}

{\theorembodyfont{\rmfamily} \theoremheaderfont{\scshape} \theoremseparator{.}
\theoremstyle{break}

\newtheorem{remarks}[theorem]{Remarks}}

\title{On finite dimensionality of mixed Tate motives}
\author{{Shahram Biglari}~\footnote{Email Address: \texttt{biglari@math.uni-bielefeld.de}}\\{\footnotesize Fakult\"at f\"ur Mathematik, Universit\"at Bielefeld, D-33501 Bielefeld, Germany}}
\date{}
\long\def\symbolfootnote[#1]#2{\begingroup%
\def\thefootnote{\fnsymbol{footnote}}\footnote[#1]{#2}\endgroup}
\begin{document}
\maketitle
\symbolfootnote[0]{\emph{2000 Mathematics Subject Classification}. Primary  11G09 - Secondary   19A99.}
\symbolfootnote[0]{\textit{Key words and phrases}. Tate motive, Grothendieck group.}
\begin{abstract}
We prove a few results concerning the notion of finite dimensionality of mixed Tate motives in the sense of Kimura and O'Sullivan. It is shown that being oddly or evenly finite dimensional is equivalent to the vanishing of certain cohomology groups defined by means of the Levine weight filtration. We then explain the relation to the Grothendieck group of the triangulated category $D$ of mixed Tate motives. This naturally gives rise to a $\lambda-$ring structure on $K_0(D)$.
\end{abstract}
\section{Introduction}\label{introduction}
The triangulated category of motives we consider is that of geometric mixed motives $DM_{\text{gm}}(\text{Spec}(k),{\mathbbm Q})$ over a perfect field $k$ constructed by Voevodsky~\cite{FSV-5}. Everything below will also work for Levine's triangulated category of mixed motives constructed in~\cite{levine-mixed}. Moreover by~\cite[Chapter VI, 2.]{levine-mixed}, the two triangulated categories above are equivalent in characteristic zero. In particular, following Levine~\cite{levine-93}, we consider the tensor triangulated full subcategory $D$ of $DM_{\text{gm}}$ generated by the Tate motives ${\mathbbm Q}(n)$, $n\in {\mathbbm Z}$.\\

This article consists of two parts. In the first part (i.e. \S\S 2-4) we prove certain results on the finite dimensionality of mixed Tate motives. In the second part (i.e. \S 5) we study the Grothendieck group of the triangulated category of mixed Tate motives with applications to problems discussed in the first part. Now we explain in more detail the content of this article.\\

The triangulated category $D$ of mixed (Tate) motives is a ${\mathbbm Q}-$linear tensor triangulated category. We recall in \S 2 the basic properties of such categories. These are enough for the following basic definitions. For example, it follows from these properties that we have a functorial and well-defined action of the symmetric group $\Sigma_n$ (i.e. the set of bijective functions on $\{1,2,\cdots,n\}$) on the $n-$th tensor power of a motive $X$. Moreover we show~(\ref{lambda-sign}) that this action is well-behaved with respect to the shift functor $X\mapsto X[1]$.\\

The definition of Schur functors can be given in any ${\mathbbm Q}-$linear tensor category in which projectors have a kernel, i.e. Karoubian categories. This is done by Deligne~\cite{deligne-2002}. We show in~(\ref{dmt-is-karoubian}) that the triangulated category of mixed Tate motive is Karoubian. This allows us to use the most basic properties~\cite[Section 1]{deligne-2002} of Schur functors $X\mapsto S_{\lambda}(X)$. the good behaviour of the action of $\Sigma_n$, as described above, gives a functorial isomorphism
\[
S_{\lambda}(X[1])\simeq S_{\lambda^t}(X)[n]
\]
where $\lambda$ is a partition of $n$ and $\lambda^t$ is the transpose partition. This is proved in detail in~(\ref{schur=translation}).\\

A more important structure of the triangulated category of mixed Tate motives is the existence of a weight filtration in the sense of Be{\u\i}linson, Ginsburg, and Schechtman~\cite[1.3.1]{bgs}. This is the result of Levine~\cite{levine-93}. Using this we explain in~(\ref{conservativeness-2}) the existence of a conservative tensor functor $\overline{\text{gr}}^W$ on $D$ with values in the bounded derived category of finite dimensional ${\mathbbm Q}-$vector spaces. Therefore we have the isomorphism
\[
S_{\lambda}\bigl(\overline{\text{gr}}^W(X)\bigr)\simeq \overline{\text{gr}}^W\bigl(S_{\lambda}(X)\bigr)
\]
of finite dimensional vector spaces. Denote by $d^+(X)$ (resp. $d^-(X)$) the sum of dimensions of even (resp. odd) cohomology ${\mathbbm Q}-$vector space of $\overline{\text{gr}}^W(X)$. We prove the following theorem.
\begin{no-theorem}[{\ref{vanishing-of-alt}}]
Let $X$ be a mixed Tate motive and $\mu$ a partition. Then $S_\mu(X)=0$ if and only if $[\mu]\supseteq [1,d^+(X)+1]\times [1,d^-(X)+1]$.
\end{no-theorem}
In particular, for a mixed Tate motive $X$ the following assertions are equivalent.
\begin{enumerate}
\item There exists $n\geq 1$ such that $\text{Alt}^n(X)=0$.
\item $H^q\bigl(\overline{\text{gr}}^W(X)\bigr)=0$ for every $q\equiv 1 (\text{mod }2)$.
\end{enumerate}
Similarly, there exists $n\geq 1$ such that $\text{Sym}^n(X)=0$ if and only if $H^q\bigl(\overline{\text{gr}}^W(X)\bigr)=0$ for every $q\equiv 2 (\text{mod }2)$.\\

We mention two consequences of this result. Recall from Kimura~\cite{kimura} that a motive $X$ is called $\text{Alt}$ (or evenly) finite dimensional if the equivalent conditions in the theorem are satisfied, i.e. if there exists an integer $n$ with $\text{Alt}^n(X)=0$. Similarly $\text{Sym}$ (oddly) finite dimensionality is defined. The above theorem implies~(\ref{2out3}) that if $X\to Y\to Z\to X[1]$ is a distinguished triangle in which $X$ and $Z$ are $\text{Sym}$ finite dimensional, then so is $Y$. Similar statement holds for $\text{Alt}$ finite dimensional objects. This corollary is also a special case of a more general result by Mazza~\cite{mazza} and Guletski{\u\i}~\cite{guletskii}. Let $X$ be a mixed Tate motive which is $\text{Alt}$ (respectively $\text{Sym}$) finite dimensional. Let us define the dimension of such a motive $X$, denoted $\text{dim}_{\pm}(X)$, to be the greatest integer $n$ such that $\text{Alt}^n(X)\neq 0$ (respectively $\text{Sym}^n(X)\neq 0$). Another consequence (\ref{vanishing-of-alt}) of the above theorem is
\[
\text{dim}_{\pm}(X)=\text{dim}_{\mathbbm Q}\bigl(H(\overline{\text{gr}}^W(X))\bigr).
\]
It is worthwhile mentioning (see below) that this integer is determined by the class $\class{}(X)$ of $X$ in $K_0(D)$.\\

A mixed Tate motive $X$ is, according to Kimura~\cite{kimura} and Andr\'e, Kahn, and O'Sullivan~\cite{andre}, said to be finite dimensional in the sense of Kimura-O'Sullivan if $X\simeq X_+\oplus X_-$ for an $\text{Alt}$ finite dimensional motive $X_+$ and a $\text{Sym}$ finite dimensional motive $X_-$. Not all mixed Tate motives are finite dimensional in this sense. In~(\ref{examplenonfinite}) we give an infinite set of non-isomorphic non-finite dimensional motives. However it follows~(\ref{vanishing-of-alt}) that for any mixed Tate motive $X$ there exists an integer $n$ such that the Schur functor $S_{n\times n}$ vanishes on $X$ where $n\times n$ is the partition $(\lambda_1,\cdots,\lambda_n)$ with $\lambda_i=n$. For a bibliography and survey of results and conjectures concerning finite dimensionality of motives see Andr\'e's article~\cite{andre0}.\\

The second part of this article is devoted to a brief study of the Grothendieck group of $D$. We show in~(\ref{K0-ring}) that $K_0(D)$, which is in fact a ring because of the existence of a bi-exact tensor structure on $D$, is as a ring isomorphic to ${\mathbbm Z}[\tau,\tau^{-1}]$ where $\tau$ is an indeterminate. An isomorphism is given by $\tau\mapsto \class{}({\mathbbm Q}(1))$. We have the following result.
\begin{no-theorem}[\ref{lambda-structure}]
For each $i\geq 0$ the map
\[
\lambda^i\colon K_0(D)\to K_0(D),\quad \class{}(X)\mapsto \class{}(\Alt{i}{X})
\]
is well-defined. The $\lambda^i$ define a structure of $\lambda-$ring on $K_0(D)$.
\end{no-theorem}
Essentially the same proof as that of~\ref{lambda-structure} gives in fact a more general statement~(\ref{general-lambda-structure}), in which arbitrary Schur functors are considered. These results imply in particular the rationality of the Zeta function $\zeta_X(t)$ of a mixed Tate motive $X$, defined following~\cite{kapranov} and~\cite{andre0}*{4.3} as an element of $K_0(D)[\![t]\!]$.\\

It remains to make a remark in the case of (not necessarily Tate) geometric motives. It would be useful and interesting to have a list of basic properties similar to ones above (particularly the previous theorem) using instead the slice filtration on motives as defined and studied by Voevodsky~\cite{voe-98} and Huber and Kahn~\cite{huber-kahn}. The basic problem is to determine the Grothendieck group of the triangulated category where graded pieces of this filtration land. This is a work in progress.\\

\emph{Acknowledgements}. The first part (i.e. sections 1-4) of this article is an expanded version of a section of author's Ph.~D. dissertation. I thank my supervisor Annette Huber-Klawitter for her helps and encouragements. I would also like to thank the IH\'ES for its support and hospitality during the period January-March 2006 where the work on the last section of this article was done. Special thanks goes to the referee for his careful comments and contributions which improved the article a lot. 
\section{Mixed Tate motives}\label{mtm}
Similar to Levine's definition~\cite[Definition 3.1]{levine-93}, define the triangulated category of mixed Tate motives $D:=DMT({\rm Spec}(k), {\mathbbm Q})$ to be the tensor triangulated full-subcategory of Voevodsky geometric motives $T:=DM_{{\rm gm}}({\rm Spec}(k), {\mathbbm Q})$ of~\cite{FSV-5} generated by ${\mathbbm Q}(\pm 1)$. The categories $D$ and $T$, being tensor triangulated, satisfy the following axioms.
\begin{list}{}{\labelwidth=50pt \leftmargin=50pt}
\item[(TTC$_\otimes$)] $D$ has the structure $(\otimes, \varphi, t, {\mathbbm 1})$ of an $ACU$(= compatible associative, commutative, and unital) ${\mathbbm Q}-$linear tensor category in the sense of~\cite{saavedra-72}*{Chapitre I, 2.4.3}.
\item[(TTC$_{\triangle}$)] $D$ has the structure of a ${\mathbbm Q}-$linear triangulated category.
\item[(TTC$_{\rm R,L}$)] There exist natural isomorphisms $\rho_{X,Y}\colon X\otimes Y[1]\to (X\otimes Y)[1]$ and $\lambda_{X,Y}\colon X[1]\otimes Y\to (X\otimes Y)[1]$ of functors from $D\times D\to D$ such that for every distinguished triangle $X\to Y\to Z\to X[1]$ with differential $d\colon Z\to X[1]$ the resulting triangles $T\otimes X\to T\otimes Y\to T\otimes Z \to (T\otimes X)[1]$ with differential $\rho_{T,X}\circ ({\rm id}_T\otimes d)$ and $X\otimes T\to Y\otimes T\to Z\otimes T \to (X\otimes T)[1]$ with differential $\lambda_{X,T}\circ (d\otimes {\rm id}_T)$ are distinguished.
\item[(TTC$_{-1}$)] the diagram
$$
\xymatrix{
\mathbbm{1}[1]\otimes\mathbbm{1}[1]\ar[d]^-{l_{{\mathbbm 1},{\mathbbm 1}}}\ar[r]^-{{\rm twist}} & \mathbbm{1}[1]\otimes\mathbbm{1}[1]\ar[d]^-{l_{{\mathbbm 1},{\mathbbm 1}}}\\
\mathbbm{1}[2]\ar[r]^-{-{\rm id}} & \mathbbm{1}[2]}
$$
commutes.
\end{list}
\begin{remark}
For a discussion on the above axioms and other notions of tensor triangulated categories see~\cite[\S 2]{shahram-2007}. For the case of triangulated category of geometric motives of Voevodsky see~\cite[2.1.3]{FSV-5} and also~\cite[Appendix 8]{mvw}.
\end{remark}
Let $D$ be a category satisfying the above axioms and $X$ an object of $D$. Existence and functoriality of the commutativity constraint gives a natural action of the symmetric group $\Sigma_2$ on $X^{\otimes 2}$. This can be uniquely extended to an action of the symmetric group $\Sigma_n$ on $X^{\otimes n}$. In particular there is an algebra homomorphism
$$
\xi_X\colon {\mathbbm Q}[\Sigma_n]\to {\rm End}_D(X^{\otimes n})
$$
which is natural in $X$, i.e. $\xi_Y(\pi)\circ f^{\otimes n}=f^{\otimes n}\circ \xi_Y(\pi)$ for any morphism $f\colon X\to Y$ and any $\pi\in {\mathbbm Q}[\Sigma_n]$. We usually denote $\xi_X(\pi)$ by $\pi_X$ or just $\pi$. We remark that this action is well-defined only because of the compatibility conditions (i.e. coherence theory) in the definition of tensor categories.
\begin{lemma}\label{lambda-sign}
There exists a functorial isomorphism $\alpha_{X,n}\colon (X[1])^{\otimes n}\to X^{\otimes n}[n]$ such that for each $\sigma\in  \Sigma_n$ the diagram
$$
\xymatrix{
(X[1])^{\otimes n}\ar[r]^-{\alpha_{X,n}}\ar[d]^-{\sigma_{X[1]}} & X^{\otimes n}[n]\ar[d]^-{\epsilon(\sigma)\sigma_{X}}\\
(X[1])^{\otimes n}\ar[r]^-{\alpha_{X,n}} & X^{\otimes n}[n]
}
$$
is commutative.
\end{lemma}
\begin{proof}
First consider the case $n=2$ and $\sigma_X =t\colon X^{\otimes 2}\to X^{\otimes 2}$. Let $u_X\colon X\to {\mathbbm 1}\otimes X$ be the natural isomorphism. For objects $Z, T$ consider the diagram
$$
\xymatrix{
Z[1]\otimes T[1]\ar[d]^-{t} \ar[rr]^-{u_{Z}[1]\otimes u_{T}[1]} & & ({\mathbbm 1}\otimes Z)[1]\otimes ({\mathbbm 1}\otimes T)[1]\ar[rr]^-{\lambda_{{\mathbbm 1}, Z}^{-1}\otimes \lambda_{{\mathbbm 1}, Z}^{-1}}\ar[d]^-{t} & & ({\mathbbm 1}[1]\otimes Z)\otimes ({\mathbbm 1}[1])\otimes T)\ar[d]^-{t}\\
T[1]\otimes Z[1] \ar[rr]^-{u_{T}[1]\otimes u_{Z}[1]} & & ({\mathbbm 1}\otimes T)[1]\otimes ({\mathbbm 1}\otimes Z)[1]\ar[rr]^-{\lambda_{{\mathbbm 1}, T}^{-1}\otimes \lambda_{{\mathbbm 1}, T}^{-1}} & & ({\mathbbm 1}[1]\otimes T)\otimes ({\mathbbm 1}[1])\otimes Z).
}
$$
This is commutative by functoriality of the commutativity constraint $t$. Next consider the diagram
$$
\xymatrix{
({\mathbbm 1}[1]\otimes Z)\otimes ({\mathbbm 1}[1]\otimes T)\ar[rr]^-{\id{}\otimes t\otimes \id{}}\ar[d]^-{t} & & ({\mathbbm 1}[1]\otimes {\mathbbm 1}[1])\otimes (Z\otimes T)\ar[d]^-{t\otimes t}\\
({\mathbbm 1}[1]\otimes T)\otimes ({\mathbbm 1}[1]\otimes Z)\ar[rr]^-{\id{}\otimes t\otimes \id{}} && ({\mathbbm 1}[1]\otimes {\mathbbm 1}[1])\otimes (T\otimes Z).
}
$$
This is a formal diagram and hence, by pentagon and hexagon axioms, is commutative in any tensor category, see also~\cite[Chapitre I, apr\`{e}s 2.5.3.4]{saavedra-72}. Next let $f\colon Y\to Y'$ be a morphism. The diagram
$$
\xymatrix{
({\mathbbm 1}[1]\otimes {\mathbbm 1}[1])\otimes Y\ar[d]^-{({\rm id}\otimes {\rm id})\otimes f} \ar[r]^-{\varphi^{-1}} & {\mathbbm 1}[1]\otimes ({\mathbbm 1}[1]\otimes Y)\ar[d]^-{{\rm id}\otimes ({\rm id}\otimes f)} \ar[r]^-{\id{}\otimes\lambda} & {\mathbbm 1}[1]\otimes Y[1]\ar[r]^-{\lambda}\ar[d]^-{{\rm id}\otimes f[1]}& Y[2]\ar[d]^-{f[2]}\\
({\mathbbm 1}[1]\otimes {\mathbbm 1}[1])\otimes Y' \ar[r]^-{\varphi^{-1}} & {\mathbbm 1}[1]\otimes ({\mathbbm 1}[1]\otimes Y') \ar[r]^-{\id{}\otimes\lambda} & {\mathbbm 1}[1]\otimes Y'[1]\ar[r]^-{\lambda}& Y'[2]
}
$$
in which $\varphi^{-1}$ is the associativity isomorphism, is commutative by functoriality of $\lambda$ and ${\varphi}$. To prove the lemma in the first case, put all these diagrams together and let $Z=T=X$, $Y=Y'=X^{\otimes 2}$, and $f=-t$. For $n>2$, note that each $\sigma$ can be written as a product of transpositions. The general diagrams above and well-definedness of $\xi_X$ show that the statement holds for any $n$ and $\sigma$.
\end{proof}

\section{Levine weight filtration}
In~\cite[\S 1]{levine-93} Levine defines a weight filtration on his triangulated category of mixed Tate motives. Below we mention this result, when imitated in the setting of Voevodsky triangulated category of geometric motives, and use it to obtain more specific properties of mixed Tate motives. These are used in next section to define, for example, the alternating algebra object of an object $X$ of $D$ inside the category $D$. In what follows $D=DMT({\rm Spec}(k),\mathbbm{Q})$ is the tensor triangulated category of mixed Tate motives with rational coefficients over a perfect field $k$ as in~\S\ref{mtm}.
\begin{definition}\label{weight-categories}
For an interval $J\subseteq \mathbbm{Z}$ define $W_{J}D$ to be the full triangulated subcategory of $D$ generated by $\mathbbm{1}(m)$ with $m\in J$.
\end{definition}
\begin{theorem}[Levine]\label{weight-filtration-start}
Let $D=DMT({\rm Spec}(k),\mathbbm{Q})$ be the tensor triangulated category of mixed Tate motives and $n$ an integer.
\begin{enumerate}[label=(\roman*), ref=\thetheorem(\roman*)]
\item\label{weight-derived} The category $W_{[n,n]}D$ is, as a triangulated category, equivalent to the bounced derived category of finite dimensional $\mathbbm{Q}-$vector spaces.
\item The inclusion $W_{]-\infty,n[}D\to D$ has a left adjoint $W_{< n}\colon
D\to W_{]-\infty, n[}D$.
\item The inclusion $W_{]m,+\infty[}D\to D$ has a right adjoint $W_{> m}\colon
D\to W_{]m,+\infty[}D$.
\item\label{weight-triangles-main} For each object $X$ the canonical triangle $W_{>n-1 }X\to X\to W_{<n}X\to W_{> n-1}X[+1]$ is distinguished and natural in $X$.
\end{enumerate}
\end{theorem}
\begin{proof}
This is proved in \cite[\S 1]{levine-93}.
\end{proof}
\begin{remarks}\label{remarks-weight}
\begin{enumerate}[label=(\arabic*), ref=\thetheorem(\arabic*)]
\item\label{triangulated-weight} We note that the functors $W_{<n}$ and $W_{>m}$ being adjoint to triangulated functors are triangulated. By abuse of notation we occasionally consider these as endofunctors on $D$.
\item $W_{<n}\circ W_{>m}=W_{>m}\circ W_{<n}$.
\item\label{weight-of-generators} If $X={\mathbbm 1}(q)[p]$ then $W_{<n}(X)$ is $X$ if $q<n$ and zero otherwise. Similarly $W_{>m}(X)$ is $X$ if $q>m$ and zero otherwise. This follows from~\ref{weight-filtration-start} and in fact used as definition of $W_{<n}(X)$ in \cite[\S 1]{levine-93}.
\item\label{weight-example} Let $S\to \text{Spec}(k)$ be a connected smooth projective equidimensional scheme of dimension $d\geq 1$. Assume that $X=M(S)$ is mixed Tate. Then it is easy to see that $W_{<0}(X)=0$, $W_{<1}(X)={\mathbbm 1}$, and $W_{>d}(X)=0$.
\end{enumerate}
\end{remarks}
\begin{definition}
${\rm gr}^W_n(X):=W_{>n-1}\circ W_{<n+1}(X)$ and $\text{gr}^W(X):=\coprod \text{gr}^W_n(X)$
\end{definition}
\begin{cor}\label{weight-triangles}
For each $n\in \mathbbm{Z}$, the triangles ${\rm gr}^W_n(X)\to W_{<n+1}(X)\to W_{<n}(X)\to {\rm gr}^W_n(X)[+1]$ and $W_{> n}(X)\to W_{> n-1}(X)\to {\rm gr}^W_n(X)\to W_{> n}(X)[+1]$ are distinguished and functorial in $X$.
\end{cor}
\begin{proof}
This follows from~\ref{weight-filtration-start} and \ref{triangulated-weight}.
\end{proof}
\begin{cor}\label{finite-gr}
Let $X$ be a mixed Tate motive. The set of integers $n$ for which ${\rm gr}^W_n(X)\neq 0$ is finite.
\end{cor}
\begin{proof}
To prove this note that by~\ref{triangulated-weight} the endo-functor ${\rm gr}^W_n$ is triangulated (i.e. exact) and hence it is enough to prove the statement for $X$ of the form $\mathbbm{1}(q)[p]$. But this follows from~\ref{weight-of-generators}.
\end{proof}
\begin{definition}\label{def-grw}
$\overline{\rm gr}^W\colon D\to D^b({\mathbbm Q}-{\rm fgmod})$ is defined by ${\displaystyle{X\mapsto \underset{n}{{\coprod}}{\rm gr}_n^W(X)(-n)}}$.
\end{definition}
\begin{proposition}\label{conservativeness-2}
$\overline{\rm gr}^W$ is a tensor triangulated functor and $\overline{\rm gr}^W(X)=0$ if and only if $X=0$.
\end{proposition}
\begin{proof}
Let $X$ be a mixed Tate motive. Since ${\rm gr}_n^W(X)=0$ for almost all $n$, the functor $\overline{\rm gr}^W$ is well-defined and is obviously additive. Let $\overline{\rm gr}^W(X)=0$. By definition ${\rm gr}_n^W(X)=0$ for all $n\in {\mathbbm Z}$. The result~\ref{finite-gr} implies that $X=0$. To see that $\overline{\rm gr}^W$ is an exact (i.e. triangulated) functor note that by~\ref{triangulated-weight} the individual functors ${\rm gr}_n^W$ are exact. Also the functor $X\mapsto X(-n)$ is triangulated. It follows that the functor $\overline{\rm gr}_n^W$ is also triangulated. Since a coproduct of distinguished triangles is distinguished, it follows that $\overline{\rm gr}^W$ is a triangulated functor. Now we prove that $\overline{\rm gr}^W$ is a tensor functor. First note that $\overline{\rm gr}^W({\mathbbm 1})$ is the complex ${\mathbbm Q}$ concentrated in degree $0$. Let $Y$ be another mixed Tate motive. We show that there is a natural isomorphism
$$
\xi_{X,Y}\colon \overline{\rm gr}^W(X)\otimes \overline{\rm gr}^W(Y)\to \overline{\rm gr}^W(X\otimes Y)
$$
of functors. To see this we first construct a morphism $\gamma_{p,q}\colon \text{gr}_p^W(X)\otimes \text{gr}_p^W(X)\to \text{gr}_n^W(X\otimes Y)$ for any $p+q=n$ and any pair of objects $X,Y$ which is natural in both arguments and which is an isomorphism whenever $Y$ is an object of $W_{[q,q]}D$. Consider the natural morphisms
$$
\xymatrix{
X\otimes Y & W_{>p-1}(X)\otimes W_{>q-1}(Y)\ar[l]_-{\alpha_{p,q}}\ar[r]^-{\beta'_{p,q}} &  \text{gr}^W_p(X)\otimes W_{> q-1}(Y)\ar[r]^-{\beta''_{p,q}} & \text{gr}^W_p(X)\otimes \text{gr}^W_q(Y).
}
$$
Using the canonical distinguished triangle $W_{> p}(X)\to W_{> p-1}(X)\to {\rm gr}^W_p(X)\to W_{> p}(X)[+1]$ in~\ref{weight-triangles}, tensoring it with $W_{> q-1}(Y)$, and noting that $W_{<n+1}(W_{> p}(X)\otimes W_{> q-1}(Y))$ vanishes we see that $\text{gr}^W_n(\beta')$ is an isomorphism. Similarly $\text{gr}^W_n(\beta'')$ is an isomorphism. Now define
$$
\gamma_{p,q}:=\text{gr}^W_n(\alpha_{p,q})\circ \text{gr}^W_n(\beta')^{-1}\circ \text{gr}^W_n(\beta'')^{-1}.
$$
By definition $\gamma_{p,q}$ is an isomorphism whenever $X=\text{gr}^W_p(X)$ or $Y=\text{gr}^W_q(Y)$. Denote the coproduct of all $\gamma_{p,q}(-p-q)$ by $\xi_{X,Y}$. This is natural in $X, Y$ and induces isomorphism for generators of $D$, we conclude that $\xi_{X,Y}$ is an isomorphism for all $X, Y$. These isomorphisms are compatible with commutativity, unit and associativity constraints. To see this first note that, as above, the fact that $W_{>n}$ is right adjoint to the inclusion $W_{>n}D\to D$ implies that the morphism $\alpha_{X, Y}$ gives a morphism
\[
\rho_{X,Y}\colon W_{\geq p}(X)\otimes W_{\geq q}(Y)\to W_{\geq n}(X\otimes Y)
\]
where $n=p+q$ and $W_{\geq n}=W_{>n-1}$. This is natural in $X$ and $Y$ and compatible with tensor structures. For example let $t$ be the commutativity constraint. We must show that the diagram
\[
\xymatrix{
W_{\geq p}(X)\otimes W_{\geq q}(Y)\ar[r]^-{\rho_{X,Y}}\ar[d]^{t} & W_{\geq n}(X\otimes Y)\ar[d]^{W_{\geq n}(t)}\\
W_{\geq q}(Y)\otimes W_{\geq p}(X)\ar[r]^-{\rho_{Y,X}} & W_{\geq n}(Y\otimes X)
}
\]
is commutative. This diagram being obtained as the right adjoint of the commutative diagram
\[
\xymatrix{
W_{\geq p}(X)\otimes W_{\geq q}(Y)\ar[r]^-{\alpha_{p,q}}\ar[d]^{t} & X\otimes Y\ar[d]^{t}\\
W_{\geq q}(Y)\otimes W_{\geq p}(X)\ar[r]^-{\alpha_{q,p}} & Y\otimes X,
}
\]
is commutative. In short, we have commutative diagrams, analogous to those for a "lax monoidal functor", which express the compatibility of the above morphism with the unit, associativity and commutativity constraints. Similarly we have a canonical morphism
\[
{\sigma_{X,Y}}:W_{\leq n}(X\otimes Y)\to W_{\leq p}(X)\otimes W_{\leq q}(Y)
\]
where $n=p+q$ and $W_{\leq n}=W_{<n+1}$. This is natural in $X$ and $Y$ and there are commutative diagrams analogous to those for a "colax monoidal functor". Then $\gamma_{p,q}$ is the composite
\[
\xymatrix{
W_{\leq n}W_{\geq n}(X\otimes Y)& \ar[l]_-{s} W_{\leq n}\bigl(W_{\geq p}(X)\otimes W_{\geq q}(Y)\bigr) & W_{\leq p}W_{\geq p}(X)\otimes W_{\leq q}W_{\geq q}(Y)\ar[l]_-{\simeq}^-{r}.
}
\]
of $r=({\sigma_{W_{\geq p}(X),W_{\geq q}(Y)}})^{-1}$ with $s=W_{\leq n}(\rho_{X,Y})$. The necessary compatibilities for the morphisms $\gamma_{p,q}$ can now be deduced from those for $\rho$ and $\sigma$ in a similar way to those for the composite of two tensor functors. The result follows.
\end{proof}
\begin{cor}\label{grf-iso}
Let $f\colon X\to Y$ be a morphism of mixed Tate motives. Then $f$ is an isomorphism if and only if $\overline{\text{gr}}^W(f)$ is.
\end{cor}
\begin{proof}
The only if part follows from functoriality of $\overline{\text{gr}}^W_n$. To prove the if part, let $Z$ be a mixed Tate motive such that $X\to Y\to Z\to X[+1]$ is a distinguished triangle. The functor $\overline{\text{gr}}^W$ takes a distinguished triangle to a distinguished one. This implies that $\overline{\text{gr}}^W(Z)=0$. The previous result shows that $Z=0$ and hence $f$ is an isomorphism.
\end{proof}
\begin{cor}\label{xy=0}
If $X\otimes Y=0$, then $X=0$ or $Y=0$.
\end{cor}
\begin{proof}
The same is true for $D^b({\mathbbm Q}-{\rm fgmod})$ and hence by~\ref{conservativeness-2} for mixed Tate motives.
\end{proof}
\begin{proposition}\label{dmt-is-karoubian}
$D=DMT(\text{Spec}(k),{\mathbbm Q})$ is Karoubian, i.e. every projector has a kernel.
\end{proposition}
\begin{proof}
Let $p\colon X\to X$ be a projector. It follows from~\ref{weight-filtration-start} that for each $n$ the projector $\text{gr}^W_n(p)$ has a kernel. An induction based on~\ref{weight-of-generators} shows that it is enough to prove that if $W_{>n-1}(p)$ and $W_{<n}(p)$ have kernel, then so does $p$. View $D$ as a subcategory of the triangulated category $DM_{\text{gm}}$ of geometric motives. The diagram
$$
\xymatrix{
W_{>n-1 }X\ar[r]\ar[d]^-{W_{>n-1}(p)} & X\ar[r]\ar[d]^-{p} &W_{<n}X\ar[r]\ar[d]^-{W_{<n}(p)} &W_{> n-1}X[+1]\ar[d]^-{W_{>n-1}(p)[1]}\\
W_{>n-1 }X\ar[r] & X\ar[r] &W_{<n}X\ar[r] &W_{> n-1}X[+1]
}
$$
is a morphism of distinguished triangles and each vertical morphism is a projector. The category $DM_{\text{gm}}$ is by definition Karoubian. Note that the triangle
$$
\text{ker}_{DM_{\text{gm}}}(W_{>n-1}(p))\to \text{ker}_{DM_{\text{gm}}}(p)\to \text{ker}_{DM_{\text{gm}}}(W_{<n}(p))\to \text{ker}_{DM_{\text{gm}}}(W_{>n-1}(p))[1]
$$
being a direct summand of the previous distinguished triangle is a distinguished triangle in $DM_{\text{gm}}$. This means that the kernel of $p$ is isomorphic to an object of $D$.
\end{proof}
\begin{cor}\label{thickness}
$D$ is a thick tensor subcategory of the triangulated category of geometric motives.
\end{cor}
\begin{proof}
Let $X,Y$ be geometric motives such that $Z:=X\oplus Y$ is mixed Tate. Note that $Y$ is the kernel of an idempotent of $Z$. The result follows from~\ref{dmt-is-karoubian}.
\end{proof}
\begin{cor}
Let $S$ be a smooth scheme of finite type over $k$ and $n\geq 1$ an integer. If $M_{gm}(S)^{\otimes n}$ is mixed Tate, then so is $\Mgm{S}$.
\end{cor}
\begin{proof}
We assume that $n\geq 2$. Let $\Delta\colon S\to S^{\times_{k}n}$ be the diagonal morphism and $\text{pr}_1\colon S^{\times_kn}\to S$ the projection onto the first factor. Note that $\Mgm{S}$ is the image of the projector $M_{\text{gm}}(\Delta)\circ M_{\text{gm}}(\text{pr}_1)$. The result follows from~\ref{thickness}.
\end{proof}
\section{Schur functors and finite dimensionality}
The construction of Weyl modules and Schur functors as done, for example, in~\cite[{\sc Lecture} 6]{fulton-harris} for finite dimensional vector spaces can be carried out in any Karoubian ${\mathbbm Q}-$linear tensor category. This is explained in~\cite{deligne-2002}. Following the latter we give the definition and first properties of the Schur functors on the triangulated category of mixed Tate motives $D$ of \S\ref{mtm}.\\

Let us denote by $A$ the group algebra ${\mathbbm Q}\Sigma_n$. An approach to the theory of representation of the symmetric group $\Sigma_n$, due to Young, is based on a statement that there exists a bijection between the set of isomorphism classes of irreducible $A-$modules and that of partitions of $n$, i.e. sequences $(\lambda_1,\ldots,\lambda_k)$ of positive integers $\lambda_1\geq \lambda_2\geq\cdots$ the sum of whose elements is $n=:|\lambda|$. See~\cite{fulton-harris} or~\cite{james} for details. For a partition $\lambda$ of an integer $n$ define its diagram to be
$$
[\lambda]:=\{(i,j)\in{\mathbbm Z}^2\ |\ 1\leq i, 1\leq j\leq \lambda_i\ \}\subseteq [1,n]\times [1,n].
$$
Denote by $V_\lambda$ the irreducible $A-$module corresponding to $\lambda$. We note that $V_\lambda$ is in fact defined over ${\mathbbm Q}$, i.e. $\text{End}_A(V_\lambda)={\mathbbm Q}$. It follows that $A\simeq \prod \text{End}_{\mathbbm Q}(V_\lambda)$, see~\cite[{\sc Lecture 4}]{fulton-harris}.\\

We recall from~\cite{deligne-2002} the following construction. Let $U_\lambda:=V_\lambda^{\vee}$ be dual of $V_\lambda$ as a ${\mathbbm Q}-$vector space. Note that $U_\lambda$ is also an $A-$module, i.e. equipped with the dual action. For an object $X$ of $D$ consider the functor
\[
F_{\lambda, X}\colon D^{\text{op}}\to {\mathbbm Q}-\text{mod}\qquad Y\mapsto \mor{\mathbbm Q}(V_\lambda, \mor{D}(Y, X^{\otimes n})).
\]
This is a representable functor. In fact we may take a ${\mathbbm Q}-$basis of $U_\lambda$ and define $U_\lambda\otimes X^{\otimes n}$ to be the direct sum of $\text{dim}_{\mathbbm Q}(V)$ copies of $X^{\otimes n}$. It follows that $U_\lambda\otimes X^{\otimes n}$ is determined up to a unique isomorphism and represents the functor $F_{\lambda, X}$. Note that this construction is functorial and hence $U_\lambda\otimes X^{\otimes n}$ has an action of $\Sigma_n$ written as $\sigma\mapsto \sigma_{U_\lambda}\otimes \sigma_X$. It is convenient to denote the representing object of $F_{\lambda, X}$ by $\underline{{\rm Hom}}_{D}(V_{\lambda}, X^{\otimes n})$ instead of $U_\lambda\otimes X^{\otimes n}$.\\

Now let $Y$ be an object of $D$ equipped with an action of $\Sigma_n$, i.e. a ${\mathbbm Q}-$algebra homomorphism $\xi_Y\colon A\to \text{End}_D(Y)$. Note that the element $s_n:={\textstyle{\frac{1}{n!}\sum\sigma}}$ of the group ring $A$ is an idempotent. Define $Y^{\Sigma_n}$ to be the image of the projector $\xi_Y(s_n)$.
\begin{definition}[{\cite{deligne-2002}}]\label{definition-of-schur-functor} $S_\lambda\colon D\to D$ is defined to be $X\mapsto \underline{{\rm Hom}}_{D}(V_{\lambda}, X^{\otimes n})^{\Sigma_n}$.
\end{definition}
\begin{proposition}\label{schur=translation}
$S_\lambda (X[1])\simeq S_{\lambda^{t}}(X)[n]$ naturally in $X$ where the partition $\lambda^t$ is the 	transpose of $\lambda$.
\end{proposition}
\begin{proof}
Let $U'$ be ${\mathbbm Q}$ with the action $\sigma_{U'}(z)=\epsilon(\sigma)z$. Note that $U'$ is canonically identified with its dual $U'^{\vee}$. By~\ref{lambda-sign} the diagram
$$
\xymatrix{
(X[1])^{\otimes n}\ar[r]^-{\alpha_{X,n}}\ar[d]^-{\sigma_{X[1]}} & (U'\otimes X^{\otimes n})[n]\ar[d]^-{(\sigma_{U'}\otimes \sigma_X)[n]}\\
(X[1])^{\otimes n}\ar[r]^-{\alpha_{X,n}} & (U'\otimes X^{\otimes n})[n]
}
$$
is commutative. Now everything should be clear as by definition~\ref{definition-of-schur-functor} we have
$$
\begin{array}{rcl}
S_\lambda (X[1]) & = & \underline{{\rm Hom}}_{D}\bigl(V_{\lambda}, (X[1])^{\otimes n}\bigr)^{\Sigma_n}\\\\
 & = & \underline{{\rm Hom}}_{D}\Bigl(V_{\lambda},(U'\otimes X^{\otimes n})[n]\Bigr)^{\Sigma_n}\\\\
 & = & \underline{{\rm Hom}}_{D}(V_{\lambda},U'\otimes X^{\otimes n})^{\Sigma_n}[n]\\\\
 & = & \underline{{\rm Hom}}_{D}(V_{\lambda}\otimes U'^{\vee}, X^{\otimes n})^{\Sigma_n}[n]\\\\
 & = & \underline{{\rm Hom}}_{D}(V_{\lambda^{t}}, X^{\otimes n})^{\Sigma_n}[n]\\\\
 & = & S_{\lambda^{t}}(X)[n].
\end{array}
$$
Naturality of the resulting isomorphism follows from that of $\alpha_{X, n}$ and the fact that all the above isomorphisms are natural in $X$.
\end{proof}
\begin{remark}
As the proof shows, the above proposition is valid for any object $X$ of a category $D$ satisfying the axioms in \S\ref{mtm}.
\end{remark}
\begin{cor}
$\text{Alt}^n(X[1])=\text{Sym}^n(X)[n]$ and $\text{Sym}^n(X[1])=\text{Alt}^n(X)[n]$ where $\text{Alt}^n=S_{(1,\ldots,1)}$ and $\text{Sym}^n=S_{(n)}$.
\end{cor}
\begin{proposition}[{\cite{deligne-2002}}]The following assertions hold.
\begin{enumerate}[label=(\arabic*), ref=\thetheorem(\arabic*)]
\item\label{ordered-partition}If $S_{\lambda}(X)=0$ then $S_{\mu}(X)=0$ for any partition $\mu$ with $[\lambda]\subset [\mu]$.
\item\label{schur-coproduct}$S_\lambda(X\oplus Y)\ ={\displaystyle{\coprod_{|\mu|+|\eta|=|\lambda|} \bigl(S_\mu(X)\otimes S_\eta(Y)\bigr)^{[\lambda:\mu, \eta]}}}$
\end{enumerate}
\end{proposition}
Most of the basic properties of functors $S_\lambda$ follow from this proposition. We also note that if $S_\lambda (X\oplus Y)=0$, then $S_{\lambda}(X)=0$. To see this it is enough to note that $S_\lambda$, being a functor, takes the retraction $X\oplus Y\to X$ to a retraction.\\

For a mixed Tate motive $X$ let
\[
H^{\rm ev}\bigl(\overline{\rm gr}^W(X)\bigr):=\coprod_{q\in{\mathbbm Z}}H^{2q}\bigl(\overline{\rm gr}^W(X)\bigr)\quad {\rm and\quad} H^{\rm odd}\bigl(\overline{\rm gr}^W(X)\bigr):=\coprod_{q\in{\mathbbm Z}}H^{2q+1}\bigl(\overline{\rm gr}^W(X)\bigr).
\]
and define $d^+(X):=\text{dim}_{\mathbbm Q}\bigl(H^{\rm ev}\bigl(\overline{\rm gr}^W(X)\bigr)\bigr)$ and $d^-(X):=\text{dim}_{\mathbbm Q}\bigl(H^{\rm odd}\bigl(\overline{\rm gr}^W(X)\bigr)\bigr)$.
\begin{theorem}\label{vanishing-of-alt}
Let $X$ be a mixed Tate motive and $\mu$ a partition. Then $S_\mu(X)=0$ if and only if $[\mu]\supseteq [1,d^+(X)+1]\times [1,d^-(X)+1]$.
\end{theorem}
\begin{proof}
Let $X$ be as in the statement. Note that $\overline{\rm gr}^W(X)$ is an object of $D^b({\mathbbm Q}-{\rm fgmod})$. Therefore there exists an isomorphism
\[
\overline{\rm gr}^W(X)\cong\coprod_{q\in {\mathbbm Z}} H^{q}\bigl(\overline{\rm gr}^W(X)\bigr)[-q].
\]
It follows from this,~\ref{lambda-sign}, and~\ref{conservativeness-2} that the functor
\[
F\colon X\mapsto H^{\rm ev}\bigl(\overline{\rm gr}^W(X)\bigr)\oplus H^{\rm odd}\bigl(\overline{\rm gr}^W(X)\bigr)
\]
is a tensor functor from the triangulated category of mixed Tate motives to that of super finite dimensional ${\mathbbm Q}-$vector spaces and takes non-zero objects to non-zero objects. Therefore $S_\mu(X)=0$ is equivalent to $S_\mu (F(X))=F(S_\mu (X))=0$. The result~\cite{deligne-2002}*{Corollaire 1.9} states that the latter is equivalent to
\[
[\mu]\not\subseteq \{(i,j)\ |\ i\leq d^+(X)\ {\rm or}\ j\leq d^-(X)\}.
\]
The result follows immediately.
\end{proof}
\begin{cor}\label{vanishing-of-alt-2}
Let $X$ be a mixed Tate motive. The following assertions are equivalent.
\begin{enumerate}
\item There exists $n\geq 1$ such that $\text{Alt}^n(X)=0$.
\item $H^q\bigl(\overline{\text{gr}}^W(X)\bigr)=0$ for every $q\equiv 1 (\text{mod }2)$.
\end{enumerate}
Similarly, there exists $n\geq 1$ such that $\text{Sym}^n(X)=0$ if and only if $H^q\bigl(\overline{\text{gr}}^W(X)\bigr)=0$ for every $q\equiv 2 (\text{mod }2)$.
\end{cor}
\begin{cor}\label{2out3}
Let $X\to Y\to Z\to X[+1]$ be a distinguished triangle of mixed Tate motives. If $\text{Sym}^n(X)$ and $\text{Sym}^m(Z)$ reduce to zero for some non-negative integers $n, m$. Then $\text{Sym}^k(Y)=0$ for some $k\geq 0$. Similar statement holds for $\text{Alt}$.
\end{cor}
\begin{proof}
This is a direct consequence of~\ref{vanishing-of-alt-2} and the fact that $X\mapsto H^q\bigl(\overline{\rm gr}^W(X)\bigr)$ is, by~\ref{triangulated-weight}, a homological functor.
\end{proof}
\begin{remark}
The above corollary holds more generally for (not necessarily Tate) mixed motives as shown in~\cite{guletskii} and~\cite{mazza}.
\end{remark}
\begin{cor}\label{altandsym}
Let $X$ be a mixed Tate motive if $\text{Alt}^n(X)=\text{Sym}^n(X)=0$ for some $n\geq 0$, then $X=0$.
\end{cor}
\begin{proof}
This is again a direct consequence of~\ref{vanishing-of-alt-2} and~\ref{conservativeness-2}.
\end{proof}
\begin{example}\label{examplenonfinite}
Let $k$ be a field and $n$ an even positive integer, e.g. $k={\mathbbm C}$ and $n=2$, such that
$$
\mor{D}({\mathbbm 1},{\mathbbm 1}(n)[n])=K_n^M(k)\otimes_{\mathbbm Z}{\mathbbm Q}\neq 0.
$$
Consider any distinguished triangle ${\mathbbm 1}\xrightarrow{f}{} {\mathbbm 1}(n)[n]\to X\to {\mathbbm 1}[1]$ where $f\neq 0$. It follows from exactness of the functors $Y\mapsto {\rm gr}_m^W(Y)$ that ${\rm gr}^W(X)={\mathbbm 1}[1]\oplus {\mathbbm 1}(n)[n]$. Note that since $f\neq 0$ and $\mor{D}({\mathbbm 1}(n)[n],{\mathbbm 1}[1])=0$ it follows that $X\not\simeq {\rm gr}^W(X)$. Let us show that
\begin{enumerate}
\item {\sl $X$ is not finite dimensional in the sense of Kimura$-$O'Sullivan.} That is $X$ can not be written as the direct sum $Y\oplus Z$ with $\text{Alt}^p(Y)=0=\text{Sym}^q(Z)$ for some integers $p,q\geq 1$. Indeed existence of such a decomposition implies that ${\rm gr}_0^W(Z)={\mathbbm 1}[1]$ and ${\rm gr}_m^W(Z)=0$ for $m\neq 0$ and ${\rm gr}_n^W(Y)={\mathbbm 1}(n)[n]$ and  ${\rm gr}_m^W(Y)=0$ for $m\neq n$. This means that $X={\mathbbm 1}[1]\oplus {\mathbbm 1}(n)[n]$ which is a contradiction.
\item {\sl $S_{\lambda}(X)=0$ if and only if $[1,2]\times [1,2]\subseteq [\lambda]$.} This follows from~\ref{vanishing-of-alt}.
\end{enumerate}
\end{example}
\section{Interpretations in the Grothendieck group}
What follows is devoted to a study of Schur functors on the triangulated category $D$ of mixed Tate motives when viewed in the Grothendieck group. We prove the existence of a $\lambda-$ring structure on $D$. For more on $\lambda-$rings we refer to [SGA 6, Exp. V] and~\cite{knutson}. This result is useful when dealing with dimension and additivity properties of Schur functors. In the sequel we denote $W_{[n,n]}(D)$ by $D_n$ for an integer $n$. We notice that there is a natural embedding $\iota_n\colon D_n\to D$ and an exact functor $\text{gr}^W_n\colon D\to D_n$. It is clear that $\text{gr}^W_n\circ \iota_n\cong \id{}$ and $\text{gr}^W_n\circ \iota_m\cong 0$ for $m\neq n$.
\begin{lemma}\label{k0-lemma}
$D$ is an essentially small category and the class $\text{cl}(X)$ in $K_0(D)$ of each mixed Tate motive $X$ decomposes as
$$
\text{cl}(X)=\sum_{n\in {\mathbbm Z}} \text{cl}\bigl(\text{gr}^W_n(X)\bigr)=\class{}\bigl(\text{gr}^W(X)\bigr).
$$
That is $\id{}=\sum \pi_n$ where $\pi_n:=K_0(\iota_n\circ\text{gr}^W_n)$.
\end{lemma}
\begin{proof}
First note that $W_{[0,0]}$ is essentially small. Using the distinguished triangles~\ref{weight-triangles} we see by induction on $n\geq 0$ that $W_{[-n, n]}$ is essentially small. Therefore $D$, being a union of these, is also essentially small. The last assertion follows from~\ref{weight-triangles} and~\ref{finite-gr}.
\end{proof}

The Grothendieck group $K_0(D)$ is in fact a ring precisely because the tensor product $\otimes\colon D\times D\to D$ is an exact functor in both arguments. In other words
\[
K_0(D)\times K_0(D)\to K_0(D),\quad \bigl(\class{}(X),\class{}(Y)\bigr)\mapsto \class{}(X\otimes Y)
\]
defines an associative commutative unital ring structure on $K_0(D)$. Note also that the exact tensor functor $\overline{\text{gr}}^W\colon D\to D^{\text{b}}({\mathbbm Q}-\text{fgmod})$, considered in~\ref{def-grw}, gives an augmentation $\epsilon\colon K_0(D)\to {\mathbbm Z}$ defined explicitly by
\[
\epsilon (\class{}(X))={\displaystyle{\sum_{p\in{\mathbbm Z}}(-1)^p\dime{\mathbbm Q}{H^p\bigl(\overline{\text{gr}}^W(X)}\bigr)}}.
\]
\begin{proposition}\label{K0-ring}
The homomorphism $h\colon {\mathbbm Z}[\tau,\tau^{-1}]\to K_0(D)$ defined by $\tau\mapsto \text{cl}\bigl({\mathbbm 1}(1)\bigr)$ where $\tau$ is an indeterminate is an isomorphism of rings.
\end{proposition}
\begin{proof}
For an integer $n$, denote $W_{[n,n]}(D)$ by $D_n$. Let us identify ${\mathbbm Z}$ with $K_0(D_0)$. This is a ring identification by the classical K\"unneth formula. Note that the ring ${\mathbbm Z}[\tau,\tau^{-1}]$ is obtained by localising the polynomial ring ${\mathbbm Z}[\tau]$ at the multiplicative subset $\{1,\tau,\tau^2,\cdots\}$. Define
\[
h\colon {\mathbbm Z}[\tau,\tau^{-1}]\to K_0(D),\quad \tau\mapsto \text{cl}\bigl({\mathbbm 1}(1)\bigr).
\]
The function $h$ is well-defined and is a ring homomorphism. By~\ref{k0-lemma} $h$ is surjective. It remains to show that $h$ is injective. With notations as in the proof of~\ref{k0-lemma} consider the group homomorphism
\[
K_0(\text{gr}^W_n)\colon K_0(D)\to K_0(D_n).
\]
This is an split epimorphism of abelian groups. It follows that $K_0(\text{gr}^W_n)\circ h$ sends $m_q\tau^q$ to $m_q$ for $q=n$ and zero otherwise. This shows that $h$ is injective. The result follows.
\end{proof}
\begin{remarks}
\begin{enumerate}
\item The isomorphism $h\colon {\mathbbm Z}[\tau,\tau^{-1}]\to K_0(D)$ is compatible with augmentation. That is if we define $\epsilon\colon {\mathbbm Z}[\tau,\tau^{-1}]\to{\mathbbm Z}$ by $\tau\mapsto 1$, then $\epsilon\circ h=\epsilon$.
\item The result~\ref{K0-ring} gives a convenient way of thinking about mixed Tate motives. Let us denote by $[S]$ the class in $K_0(D)$ of the motive, assuming to be mixed Tate, of a smooth scheme of finite type $S\to \text{Spec}(k)$. We have the following formul\ae
$$
{[{\mathbbm A}^n_k]}=1,\quad {[{\mathbbm P}^n_k]}=1+\tau +\cdots +\tau^n,\quad{[{\mathbbm A}_k^n\setminus 0]}=1-\tau^n.
$$
\end{enumerate}
\end{remarks}
\begin{theorem}\label{lambda-structure}
For each $i\geq 0$ the map
\[
\lambda^i\colon K_0(D)\to K_0(D),\quad \class{}(X)\mapsto \class{}(\Alt{i}{X})
\]
is well-defined. The $\lambda^i$ define a structure of $\lambda-$ring on $K_0(D)$.
\end{theorem}
\begin{proof}
Using the result~\ref{K0-ring} together with~\cite{knutson}*{p. 17}, it is enough to show that the $\lambda^i$ are well-defined and define a structure of pre-$\lambda-$ring on $K_0(D)$. Let $X$ be an element of the set $\text{Obj}(D)$ of isomorphism classes of mixed Tate motives. Define
\[
\lambda (X):=\sum_{q\in{\mathbbm Z}_{+}} \class{}\bigl(\Alt{q}{X}\bigr)t^q\in 1+tK_0(D)[\![t]\!]\subset K_0(D)[\![t]\!]^{\times}.
\]
We must show that $\lambda\colon K_0(D)\to 1+tK_0(D)[\![t]\!]$ given by $\class{}({X})\mapsto \lambda (X)$ is a well-defined homomorphism of abelian groups. By definition the function $\lambda$ extends to the free abelian group on the set of isomorphism classes of mixed Tate motives, i.e.
\[
\lambda\colon {\mathbbm Z}[\text{Obj}(D)]\to 1+tK_0(D)[\![t]\!].
\]
First note that for mixed Tate motives $X$ and $Y$ we have $\lambda (\clas{}(X\oplus Y))=\lambda (\clas{}(X)+\clas{}(Y))$ where $\clas{}\colon \text{Obj}(D)\to {\mathbbm Z}[\text{Obj}(D)]$ is the canonical function. Indeed it is enough to prove that for any integer $n\geq 0$
\[
\class{}\bigl(\Alt{n}{X\oplus Y}\bigr)=\sum_{p+q=n}\class{}\bigl(\Alt{p}{X}\otimes \Alt{q}{Y}\bigr).
\]
But in view of the definition this follows from~\ref{schur-coproduct}. In particular $\lambda (\clas{}(X\oplus Y))=\lambda (X)\lambda (Y)$. Secondly we note that by~\ref{k0-lemma} we have
\[
\class{}\bigl(\Alt{p}{X}\bigr)=\class{}\bigl(\text{gr}^W(\Alt{p}{X})\bigr)=\class{}\bigl(\Alt{p}{\text{gr}^W(X)}\bigr)
\]
which implies that $\lambda(X)=\lambda(\text{gr}^W(X))$. Thirdly we show that $\lambda(X[1])=\lambda(X)^{-1}$. The previous two steps show that
\[
\lambda(X[1])=\prod_{n\in{\mathbbm Z}} \lambda(\text{gr}_n^W(X)[1]).
\]
We show that $\lambda(\text{gr}_n^W(X)[1])=\lambda(\text{gr}_n^W(X))^{-1}$ for each $n$. Using the simple fact that $\Alt{p}{X(1)}\simeq \Alt{p}{X}(p)$,~\ref{schur=translation}, and the second step we may assume that $X={\mathbbm 1}$. But the assertion is trivial in this case as 
\[
\lambda({\mathbbm 1})=1+t,\quad \lambda({\mathbbm 1}[1])=1-t+t^2-\cdots.
\]
Fourthly let $Y\to X\to Z\to Y[1]$ be a distinguished triangle. If this is split (i.e. one of the morphisms is zero), then $\lambda(X)=\lambda(Y)\lambda(Z)$ because of the first and third statements. The general case is a combination of the latter and the second step. More precisely
\[
\begin{array}{rcl}
\lambda (X) & = & \lambda\bigl(\text{gr}^W(X)\bigr)\\\\
            & = & {\displaystyle{\prod_{n\in{\mathbbm Z}}\lambda\bigl(\text{gr}_n^W(X)\bigr)}}\\\\
	    & = & {\displaystyle{\prod_{n\in{\mathbbm Z}}\Bigl(\lambda\bigl(\text{gr}_n^W(Y)\bigr)\lambda\bigl(\text{gr}_n^W(Z)\bigr)\Bigr)}}\\\\
	    & = & \lambda(Y)\lambda(Z).
\end{array}
\]
This means that $\lambda$ is a well-defined homomorphism on $K_0(D)$. 
\end{proof}
\begin{cor}[Product formula]\label{product-formula}$\lambda (x)={\displaystyle{\prod_{n\in {\mathbbm Z}}\lambda \bigl(\pi_n(x)\bigr)}}$.
\end{cor}
\begin{proof}
This was shows in the proof of~\ref{lambda-structure} and can also be derived from it using~\ref{k0-lemma} that $\id{}=\sum\pi_n$ as endomorphism of the abelian group $K_0(D)$.
\end{proof}
\begin{cor}\label{additivity-triangle-k0}
Let $Y\to X\to Z\to Y[1] $ be a distinguished triangle. Then for each integer $n\geq 0$
\[
\class{}\bigl(\Alt{n}{X}\bigr)-\class{}\bigl(\Alt{n}{Y\oplus Z}\bigr)=0\text{ in } K_0(D).
\]
\end{cor}
\begin{cor}\label{alt-sym-inverse}
${\displaystyle{\sum_{p+q=n}(-1)^q\class{}\Bigl(\Alt{p}{X}\otimes\Sym{q}{X}\Bigr)=0}}$ if $n\neq 0$.
\end{cor}
\begin{proof}
This follows from~\ref{lambda-structure} and the equation $\lambda(X\oplus X[1])=1$.
\end{proof}

By~\ref{alt-sym-inverse} we have
\[
\sum_{n\geq 0}\class{}\bigl(\Sym{n}{X}\bigr)t^n=\Bigl(\sum_{n\geq 0}\class{}\bigl(\Alt{n}{X}\bigr)(-t)^n\Bigr)^{-1}
\]
in $1+tK_0(D)[\![t]\!]$. Following~\cite{kapranov} and~\cite{andre0}*{4.3}, this element of $K_0(D)[\![t]\!]$ may be called the Zeta function of $X$ and denoted by $\zeta_X(t)$. For our purpose, let us call an element $f\in 1+tK_0(D)[\![t]\!]$ rational if there are polynomials $u$ and $v$ in $K_0(D)[t]\subseteq K_0(D)[\![t]\!]$ with $v(0)=1$ such that $f=uv^{-1}$.
\begin{cor}\label{rationality}
Let $X$ be a mixed Tate motive. Then $\zeta_X(t)$ is a rational function.
\end{cor}
\begin{proof}
The assertion is clear for objects of the form ${\mathbbm 1}(q)[p]$. Now for a distinguished triangle $X'\to X\to X''\to X[1]$ of mixed Tate motives we have $\class{}(X)=\class{}(X')+\class{}(X'')$ in $K_0(D)$. Therefore by~\ref{lambda-structure} we have the identity
\[
\zeta_X(t)=\zeta_{X'}(t)\zeta_{X''}(t).
\]
The result follows immediately.
\end{proof}
\begin{remark}
Analyzing the proof of the theorem~\ref{lambda-structure} we have the following. Let $D_{\text{ad}}$ be the underlying additive subcategory of $D$. Let $T$ be a subgroup of $K_0(D_{\text{ad}})$ generated by elements of the form
\[
\class{}\bigl(X\oplus \text{gr}^W(X)[1]\bigr)+\class{}\bigl(Y\oplus Y[1]\bigr)
\]
where $X$ and $Y$ are objects of $D_{\text{ad}}$. Then the canonical epimorphism $K_0(D_{\text{ad}})\to K_0(D)$ induces an isomorphism $K_0(D_{\text{ad}})/T\simeq K_0(D)$.
\end{remark}
It is possible to prove a slightly more general statement than~\ref{lambda-structure}. Let us first recall the folowing notations.
\begin{notations}\label{representation-ring}
For $n\geq 0$ denote by $R_n$ the Grothendieck group of the semisimple abelian category of finitely generated ${\mathbbm Q}\Sigma_n-$modules. This is a free abelian group on the set of isomorphism classes of irreducible objects or equivalently on the set of partitions of $n$. Here we adopt the convention that $0$ has one and only one partition. Define
\[
R:=\coprod_{n\geq 0} R_n.
\]
This is a commutative, associative, and unital ring with the multiplication of classes of a ${\mathbbm Q}\Sigma_p-$module $V$ and a ${\mathbbm Q}\Sigma_q-$module $W$ given by $\class{}\bigl(\text{Ind}^{\Sigma_{p+q}}_{\Sigma_{p}\times \Sigma_{q}}\bigl(V\otimes W\bigr)\bigr)$.
\end{notations}
\begin{definition}\label{def-lambda-general}
Let $X$ be an element of the set $\text{Obj}(D)$ of isomorphism classes of mixed Tate motives. Define
\[
\lambda_{\Sigma} (X):=\sum_{\mu} \class{}\bigl(S_\mu(X)\bigr)\otimes \class{}(V_\mu)t^{|\mu|}\in 1+tK_0(D)_R[\![t]\!]\subset K_0(D)_R[\![t]\!]^{\times}.
\]
where the sum is taken over all partitions $\mu$ with $|\mu|\geq 0$ and where $V_\mu$ is the irreducible representation of $\Sigma_{|\mu|}$ corresponding to $\mu$ and $K_0(D)_R:=K_0(D)\otimes_{\mathbbm Z}R$. .
\end{definition}
\begin{theorem}\label{general-lambda-structure}
$\lambda_{\Sigma}\colon K_0(D)\to 1+tK_0(D)_R[\![t]\!]$ given by $\class{}({X})\mapsto \lambda_{\Sigma} (X)$ is a well-defined monomorphism of abelian groups.
\end{theorem}
\begin{proof}
The proof is essentially the same as that of~\ref{lambda-structure}.
\end{proof}
\begin{remark}\label{lambda-remarks}
As in the special case~\ref{lambda-structure}, the previous theorem implies results similar to~\ref{product-formula} and~\ref{additivity-triangle-k0}.
\end{remark}

\end{document}